\magnification 1200

\font\tenmsb=msbm10
\font\sevenmsb=msbm10 at 7pt
\font\fivemsb=msbm10 at 5pt
\newfam\msbfam
\textfont\msbfam=\tenmsb
\scriptfont\msbfam=\sevenmsb
\scriptscriptfont\msbfam=\fivemsb
\def\Bbb#1{{\fam\msbfam\relax#1}}

 \def\la{\longrightarrow}
 \def\ni{\noindent}
 \def\cl{\centerline}
 
\def\d{\delta }
 \def\e{\epsilon}

\def\a{\alpha}
\def\b{\beta}

 \def\Y{{\cal Y}}
 \def\X{{\cal X}}
\def\W{{\cal W}}

\def\F{{\Bbb F}}
\def\N{{\Bbb N}}
\def\P{{\Bbb P}}
\def\Q{{\Bbb Q}}

 \def\Pic{{\rm Pic}}
\def\NS{{\rm NS}}

\def\Pr{{\Bbb P} ^r}
\def\Pl{{\Bbb P} ^2}

\def\Fn{{\Bbb F}_n}

\def\Vg{V_g(D)}
\def\VR{V_0(D)}
\def\Vi{V^i_0(D)}

\def\GS{V^{d,\delta }[\alpha ,\beta ]}
\def\NS{N^{d,\delta }[\alpha ,\beta ]}
\def\rGS{r^{d,\delta }[\alpha ,\beta ]}

\def\rg{r_g(D)}

\def\Ng{N_g(D)}
\def\NR{N_0(D)}
\def\Nd{N_g(d)}

\

\cl{\bf COUNTING CURVES ON SURFACES:}
\smallskip
\cl{\bf A GUIDE TO NEW TECHNIQUES AND RESULTS}

\

\

\cl { Lucia Caporaso}

\

\

\cl{ 1. INTRODUCTION}

\

\

\ni {\bf 1.1.\  Abstract and summary. \  } A series of recent results solving
classical enumerative problems for curves
on rational surfaces is described.
Impulse to the subject came from recent ideas from quantum field theory leading
to the definition of quantum cohomology. As a by-product, formulas
enumerating rational curves on certain
varieties were derived from the properties of certain generating functions
representing the free energy of
certain topological field theories. A mathematically acceptable
construction of quantum cohomology came
soon afterwards ([RT] and [KM]-[K]).
Different proofs of some of these formulas were provided later
(in [CH1] and [CH2]) using different methods that could be generalized to
cases (such as
Hirzebruch surfaces)  for which the quantum cohomology theory did not give
enumerative results.
For higher genera, the connection between
enumerative geometry  and quantum cohomology or quantum field theory is
still  largely conjectural.
 On the other hand, a recursive formula enumerating plane curves
of any genus has been recently proved using purely algebro-geometric
techniques ([CH3]). Moreover a generating
function exists together with a differential equation implying such a
recursion ([G]).

The  enumerative problem is   precisely stated in the introduction.

The second chapter contains a short description of the relation with
Quantum Cohomology and
 Kontsevich's formula for plane rational curves. A  discussion of
enumerative problems for Hirzebruch surfaces
concludes this part, which is entirely dedicated to rational curves.

The general case of plane curves of any genus is described in the third
chapter, focussing on the results of
[CH3].  The main recursive formula of that paper is explained, together with
an outline of the proof and a description of the generating
function found in [G]. At the end, there is a discussion of generalizations
of the procedure of [CH3] to
other varieties.

Chapters 2 and 3 are completely independent from each other.

In the fourth and
last chapter various
enumerative techniques  are applied to the concrete example of counting
rational plane cubics through $8$
general points.

\

\ni {\bf 1.2.\  Terminology. \  } We work with complex, projective
algebraic varieties.
   Let $S$ be
a smooth, minimal rational surface (that is,
$S$ is either  the projective plane $ \Pl$ or a  Hirzebruch  surface $\Fn$)
and let $D$ be a curve on
$S$. The linear system
$|D|$ of all curves linearly equivalent to $D$ is a projective space whose
dimension is
given by Riemann-Roch theorem (if $S$ is the plane then  $|D|$ is the
${d(d+3) \over 2}$-dimensional
projective  space of all curves of
degree
$d=\deg D$).
Moreover the genus  $p_a(D)$ of a smooth curve in
$|D|$ is constant, computable by the adjunction formula (equal to ${d-1\choose
2}$ for a smooth plane curve of degree $d$).

We  study the geometry of certain closed subvarieties of $|D|$, the
so-called ``Severi
varieties". The following notation is used throughout the paper.
Severi varieties are denoted by the symbol ``$V$", suitably decorated; for
example, by $\Vg$ we
denote the Severi variety defined to be the closure of the locus of all
irreducible curves in $|D|$
having fixed geometric genus $g$ (see below).
All of these Severi varieties have pure dimension (that is, all of their
irreducible components
have the same dimension); we  use a decorated ``$r$" for the dimension of
$V$ (equally decorated);
for example, $\rg :=\dim \Vg$.
Analogously, we use a decorated ``$N$" for the degree of $V$ as a
subvariety of $|D|$, for
example,
$\Ng := \deg \Vg$.

\

\ni {\bf 1.3.\  The problem. \  }
For a given integer $g$
(with $0\leq g\leq p_a(D)$ to avoid trivial cases),   consider  $\Vg\subset |D|$
defined above.  Also,
we denote by $V_g(d)$ the variety of irreducible plane curves of genus $g$
and degree $d$.
These varieties were first introduced
by Severi for $\Pl$ (in [S]) and have been object of much study.
For example,
it is known that their general
point represents a curve with $\d =p_a(D)-g$  nodes
(and no other singularities, a fortiori),
 and that their dimension  satisfies
 $\rg =-(K_S\cdot D) + g -1= \dim |D|-\d $, where $K_S$ is the canonical class
of $S$.

This is to say that for any $\rg $ general points on the surface $S$ there
is a finite number of
 curves of $|D|$ having geometric genus $g$ and passing through such
points. Such a number is the degree of $\Vg$ as a subvariety of $|D| $ and it is
precisely what we want to compute. More generally, we will
consider a larger class of  Severi
varieties and we will describe various ways to arrive at formulas for their
degrees.

The simplest example is that of the plane, where we are asking:

\proclaim  Question.
 How many plane curves of degree $d$ and  genus $g$ pass through
$3d+g-1$ general points?

 Until recently there was an explicit formula
answering  this
classical problem,  only in some special cases, namely
 for plane curves having few nodes (up to $6$ nodes, as
far as the author knows; due to [KP] and [Va]; see [DI] for a list of
formulas).
For example, it is not hard
to deal with  curves with only one node, and to show that for each $d$
there are $3(d-1)^2$ of them through the  appropriate number of points (see
4.1).

An interesting
recursive procedure to compute
$\Nd$ in general was suggested in $[R1]$, (but notice that the formula
there is not correct, see also [R2]
and [Ch]). Ran constructs a family where the plane degenerates to a
reducible surface (called a ``fan").
Correspondingly he gets a family of Severi varieties of which he studies
the flat limit.
This procedure  can be viewed as a recursion because such a degeneration of
$\Pl$ to a reducible surface
induces a family of irreducible curves specializing to a reducible one.

In fact, the common part to most ways of approaching these  problems
 is the use of a  technique where
the irreducible curves that one is trying to enumerate degenerate to
reducible curves,
so that one obtains an inductive formula.

\

\

\cl{2. RATIONAL CURVES}

\

\

\ni {\bf 2.1.\  Relation with Quantum Cohomology. \  } Interest and
enthusiasm for
these problems was revived by recent ideas from  quantum field theory which
led to various enumerative predictions for rational curves on varieties.

It was
proposed  by Gromov to study a series of new invariants of a given variety V.
This was also done by Witten (cf. [W]) on the base of physical intuition,
 using  intersection theory on
${\overline {M_{g,n}}}$,  the
 moduli space of Deligne-Mumford stable curves of genus
$g$, with $n$ marked points, a space well known by string theorists.

These
invariants, which  we  call 
``Gromov-Witten invariants" (they are also called ``topological $\sigma
$-models"
or ``mixed invariants" depending on context) depend
on the geometry of curves lying on
$V$.  For certain varieties (such as projective spaces) a
subclass of them  corresponds to enumerative invariants; for example, the
degrees
$N_0(d)$ for plane rational curves are   Gromov-Witten invariants.

We restrict this brief description to the case of curves of genus zero, to
avoid parts of the theory that are
still at a conjectural state, and because this is where there is a clear
link with algebraic geometry: the
above mentioned enumerative predictions all have to do with rational curves.
(See [G] for some very recent developements for genus one.)

It was conjectured (and it is now  proved)
that the Gromov-Witten invariants  satisfy a series of
properties;
 the most important of them
is  the so-called ``splitting principle" or ``composition law",
which gives a way of computing these invariants recursively.
On their existence one can
base the construction of  a family of quantum ring structures on the
cohomology ring of
$V$, deforming the standard cup product; the  associativity of this quantum
product is  a consequence of the  splitting principle.

\

\ni {\bf 2.2.\  A formula for plane rational curves. \  }
In 1993 Kontsevich derived
 a beautiful formula for rational curves in the plane,
assuming the associativity of the quantum product for $\Pl$ (not yet
proved at the time).

\proclaim Kontsevich's formula.
For $d\geq 2$
$$  N_0(d) =
\sum _{d_1 + d_2 =d} N_0(d_1)N_0(d_2)d_1d_2\left[  { 3d-4 \choose 3d_1-2
}d_1d_2 -  {
3d-4 \choose  3d_1-3}
  d_2^2 \right] .
$$

This, together with the basic fact $N_0(1)=1$, (i.e. there is a unique line
through two distinct points)
 allows
one to compute  degrees of all Severi varieties of
rational curves in the plane.

Here is, briefly, how such a formula was discovered.
Using the Gromov-Witten invariants, one defines  a generating function (the
``potential") on the cohomology ring of $V$, which completely encodes the
quantum ring structure. This is part
of the basic set-up of Quantum Cohomology, I will not say much about it and
refer to [KM] for the details.
For the special case of  $\Pl$
such a generating function is, for $\Delta \in H^*(\Pl )$:
$$
\Phi (\Delta )=\Phi ^{cl} +\Phi ^q ={1\over 2}(x_0^2x_2+x_0x_1^2) +
\sum _{d=1 }^{\infty}N_0(d) {x_2^{3d-1} \over (3d-1)! } e^{dx_1}
$$
where the variables $x_0,x_1,x_2$ are the coefficients of
$\Delta =x_0\Delta _0 +x_1\Delta _1 +x_2\Delta _2 $, with $\Delta _0$ the
identity, and
$\Delta _1$ and $\Delta _2$ the duals of the class of the line and the
point respectively.

For each $\Delta $ one gets a quantum ring structure on $H^*(\Pl )$ which
is defined using the rank-3 tensor
of all derivatives
$(\partial _j\partial _i\partial _k\Phi )_{|\Delta }$. The summand
$\Phi ^{cl}  ={1\over 2}(x_0^2x_2+x_0x_1^2)$ gives the ``generating
function" for
the classical cup product, which of course does not depend on $\Delta$

It is worth mentioning that a ``full" potential should be defined as to
include a term for  every genus
(so that $\Phi ^q$ above corresponds to the genus zero part), but it is
still an open question how to do that
in a geometrically meaningful way.

The crucial observation  is that the quantum product is associative if and
only if the potential $\Phi$
satisfies the WDVV differential equation, that is to say, if and only if
the following identity holds:
$$
\Phi ^q_{222}=(\Phi ^q_{112})^2-\Phi ^q _{111}\Phi ^q_{122}.
$$
Finally, as the reader can check by a straightforward computation,  the
above identity implies
Kontsevich's formula.

Complete proofs of Kontsevich's formula were given independently by Ruan
and Tian
(in [RT]) and by Kontsevich and Manin (in
[KM] and [K]), using rather sophisticated techniques. In both cases the
goal was to
give a mathematically rigorous definition of the Gromov-Witten invariants,
so that
they satisfy the required properties (especially the composition law!). In
[RT] this is done using symplectic
topology and the Gromov theory of pseudo-holomorphic curves. In [KM]  the authors follow an
algebro-geometric approach and use the existence of a good compactification
of the
moduli space of maps from $\P ^1 $ to $V$ (constucted in the later paper [K]).
We notice that these techniques work for a larger class of varieties; in
the specific
case of surfaces they give enumerative results for the plane, for
$\F _0 =\P ^1 \times
\P ^1$, for  $\F _1$ and for blow-ups of $\Pl$ at general points.

\

\ni {\bf 2.3.\  Hirzebruch surfaces}
In the First Reconstruction Theorem (in [KM]) there is the sketch for a
heuristic argument to obtain the
formula above.
More recently we showed  (in [CH1]) that such a  heuristic
argument
 could be made into a
completely rigorous proof involving only  classical tools, with the
advantage that this old fashioned approach
(which we call  the ``cross ratio" method) leads to formulas for rational
curves on
any rational surface $S$, that could not be found otherwise.

The methods of
[KM] and [RT] do not answer our enumerative questions for the Hirzebruch
surfaces $\Fn$. In fact
$\Fn$ is not ``convex"
  which implies that there
does not  exist a well-behaved compact moduli space of maps as it exists for
$\Pl$
(convexity here is defined follows: a variety $V$ is convex if for every
stable map $f:\P
^1 \to V$ we have $H^1(f^*T_V)=0$ where $T_V$ is the tangent bundle of $V$).
 On the other hand, the techniques  of [RT]  only depend on the symplectic
type of the
surface;  the degrees $\NR$ coincide with certain Gromov-Witten numbers as
long as $S$
is $\Pl$, $\P ^1 \times \P ^1$ or $\F _1$.
But while the Gromov-Witten numbers are symplectic-invariants, the Severi
degrees $N_0(D)$
are not, as we can see in the following example.
There are only two symplectic types of Hirzebruch
surfaces
$\Fn$, depending on the parity of $n$. So that
we can  compare
$\F _0$ with $\F _2$, which must have the same Gromov-Witten invariants.
On both surfaces, consider $D$ equal to the anticanonical
class; then one can show (using any of the techniques described in 4.1 or 4.2)
that
$N_0(-K_{\F _0}) = 12$ while $N_0(-K_{\F _2}) = 10$.

 We describe the cross ratio method in an example later (cf. 4.2).
Here we just give a summary of the results that it gives.
First, for $\Fn$, with $n\leq2$,
 and blow-ups of the plane at general points one obtains an inductive formula
(completely analogous to Kontsevich's formula)
expressing $\NR$ as a function
of simply
$N_0(D')$, with $D'<D$ (that is $D-D'$  effective and non-zero).
This corresponds to the fact that for such surfaces one can construct
degenerations of rational curves in
$|D|$ whose degenerate fibers will be reducible rational curves all of
whose components are general points of
Severi varieties $V_0(D')$.

For $n\geq 3$ a new phenomenon complicates things: this is the occurrence in
codimension one of degenerate loci that are no longer of type $V_0(D')$;
these will
instead be loci of curves satisfying certain tangency conditions. More
precisely, for
$\Fn$ let
us define the {\it tangential Severi varieties}
$\Vi \subset \VR$ as the closure of the set of curves in $\VR$ that have a
point of
contact of order $i$ with the exceptional curve $E$ of $\Fn$ (that is,
the unique curve $E$ having self intersection
$-n$). Then the cross ratio method gives a formula for $\NR$ in terms of
the degree
$N_0^i(D')$ of $V_0^i(D')$, for suitable $i$ (see [CH1]).
Therefore to have a complete picture one should also compute the degrees of
these
tangential Severi varieties.

There is another ad-hoc technique to deal with  rational curves,
 the so-called ``rational fibration method" (illustrated in example 4.2).
This is described in [CH2] where it
is applied to obtain a complete set of formulas for $\F _3$.
The picture appears to be essentially the same as for the cross ratio
approach: we
get inductive formulas expressing $\NR$ in terms of degrees of tangential Severi
varieties.
 Such a method  again focuses on the
study of one-parameter families of rational curves, and it is based on the
very basic
fact that two line bundles on
$\P ^1$ are isomorphic if they have the same degree.
See 4.2 for an example illustrating this method.

\

\

\cl{ 3. HIGHER GENERA}

\

\

\ni {\bf 3.1.\  The degeneration method of [CH3]. \  } We now consider
curves of any genus, for
which none of the methods described in the previous chapter seem to work so far.
We shall give an answer to the Question in $1.2$; from now on, we take $S=\Pl$.
Fix a line $L\subset \Pl$ once and for all, the inductive technique now
uses degenerations whose
special fiber is forced to contain an increasing number of points of $L$,
until it
becomes reducible, having to contain  $L$ itself as a component. What
remains of the
special fiber is a curve of degree $d-1$   (which may very well be reducible)
  so that one
can then use induction.

This procedure is different from the ones used for rational curves for a
crucial reason. For curves of genus
$0$ one uses generic degenerations
(such as: the family of all rational plane curves of degree $d$ through
$3d-1$ general points). Here we use a
special type of degeneration, by using the device of placing some of the
points on a fixed line.

Therefore we start by considering
a larger  class of Severi varieties,
 including certain tangential loci.
This is inspired by the case of Hirzebruch surfaces described above,  the role
of the exceptional curve $E$ is  played by the line $L$. For
rational curves on $\Fn$ we were forced to consider curves satisfying tangency
conditions with respect to $E$, getting more complicated sets of recursions.
For curves
of any genus in the plane we actually  choose to take tangential loci from the
beginning,  and this way
we get a rather simple formula.

\

\ni {\bf 3.2.\  A formula for plane  curves of any genus.\  } Let us define
generalized Severi varieties.
Let $\a = (\a_1,\ldots ,\a_h)$ and $\b = (\b_1,\ldots ,\b_k)$
be strings of nonnegative integers.
Fix  $\sum \a _j$  general points on $L\  $ denoted by
$
 \{ p_j^{(i)}\} _{1\le j \le \a_i} .
$
Assume  that $\sum i\alpha _i + \sum i\b _i =d$.
The  generalized Severi variety
$\GS$ is defined as the closure of the locus of
reduced plane curves of degree
$d$ with $\d $ nodes (hence the  geometric genus is $g = {d-1 \choose 2} -
\d$)  which

\ni
(i) do not contain $L$,

\ni (ii) have
 contact of order $i$ with $L$ at $p_j^{(i)}$ for $1\leq j \leq \a_i$
(briefly: have $\a _i$ assigned points of contact of order $i$ with $L$)

\ni
(iii) have $\b_i$ points of
contact of order $i$ with $L$ at some ``unassigned"  points.

Notice that we do not
assume the curves to be irreducible, this is why we
change notation and label by the number
of nodes
$\d$ instead of the geometric genus.

Some examples: $V^{3,1}[ (0), (3)]$
is the variety of rational cubics,
denoted by $V_0(3)$ with the notation of the previous sections
(we omit the $0$ entries in $\a$ and $\b$ whenever
this does not create confusion); and $ V^{d,\d }[(0),(d)] =V_g(d)$; also
$V^{4,3}[ (0), (4)]$ the $11$-dimensional variety of
quartics with $3$ nodes; this has two irreducible components $V_1$ and
$V_2$, where
$V_1=V_0(4)$ and $V_2$ parametrizes
all reducible quartics made of the union of a line and
a cubic.

One can show that

\ni
(1)
each irreducible component of $\GS$ has the expected dimension
$$
\rGS =\dim \GS = {d(d+3)\over 2} - \d - \sum i\a _i - \sum (i-1)\b _i =
2d+g-1 + |\b |;
$$

\ni
(2) the general point of $\GS$ parametrizes a curve having only nodes as
singularities and smooth along $L$.

Set now $\NS := \deg \GS$. Let us introduce the notation
$\a !=\a _1!\a _2 ! ...$,
${\a \choose \a'}:={\a _1\choose \a'_1}\cdot {\a _2\choose \a'_2}\cdot \ldots $
and if $S=\{ s_1, s_2, ....\}$ is any ordered set
(for example,  $S=\N$ the set of positive integers)
$S^{\a }:=s_1^{\a _1}s_2^{\a _2 } s_3^{\a _3 } \cdot \ldots$; let also $\e
^{(j)}$ be defined as
the string of integers having $1$ at the $j$th  place and $0$ elsewhere.
We can then state the main enumerative
result (See Theorem $1.1$ in [CH3]):

\proclaim Theorem.
$$
\eqalign{\NS \; &= \; \sum_{j : \b_j > 0} j
N^{d,\delta }[\alpha +\e ^{(j)},\beta -\e ^{(j)}]
\cr &\quad + \sum  \N ^{\b'-\b}{\a \choose \a'}{\b' \choose
\b}N^{d-1,\d'}[\a',\b'] \cr}
$$
where the second sum is taken over all $\a', \b'$ and $\d' \ge 0$ satisfying
$
\a' \le \a $, $
\b' \ge \b $,
$\d' \le \d $ and
$\d-\d' + |\b'-\b| = d-1 $.

\ni
The proof of this theorem uses techniques
of deformation theory together with  semistable
reduction. The degeneration technique is as follows. Let $V=\GS$; if $p\in
\Pl$ is a
point we denote by $H_p $ the hyperplane of $|D|$ parametrising curves
through $p$.
Let now $p_1,\ldots ,p_t$ be points on $L$. We consider the scheme theoretic
intersection
$$V_t :=V\cap (\cap _{i=1}^t H_{p_i})$$
which has the same degree as $V$. The formula above can be read  as a
statement describing the hyperplane
section $V\cap H_{p_1}$ as a scheme. Notice in fact that $V^{d,\delta
}[\alpha +\e ^{(j)},\beta -\e ^{(j)}] $
has codimension $1$ in $V$, while
 the  last condition $\d-\d' + |\b'-\b| = d-1 $ is precisely saying that the
codimension of $V^{d-1 ,\d'}[\a',\b']$ in $V$ must be  $1$.
The coefficients of the formula have different meanings. We will illustrate
the procedure on an example (cf. 4.3). All details can be found in [CH3].

We conclude with a nice and hopefully inspiring new way of writing the
above formula; this was found by
Getzler and appears in [G].
Let $z$  be a variable and let $u=(u_1,u_2,....)$ and $v=(v_1,v_2,....)$ be
sets of variables.
Then we define a generating function using  the degrees $\NS$
$$
G=\sum {u^{\a }\over \a !}v^{\b }\NS {z ^{\rGS } \over \rGS !}
$$

\ni
then an easy computation shows that the recursion in the above Theorem is
equivalent to the following identity

\

$$
{\partial G \over \partial z}=\sum _{k=0}^{\infty}kv_k{\partial G \over
\partial u_k}+
{\rm{Res}}_{t=0}{\rm{exp}}\bigl( \sum _{k=0}^{\infty} {u_k\over t^k} + \sum
_{k=0}^{\infty}kt^k{\partial  \over
\partial u_k}\bigr) G
$$

\

\ni
where clearly the  first summand $\sum _{k=0}^{\infty}kv_k{\partial G \over
\partial u_k}$ corresponds to the
first summand in the above formula and the remaining part corresponds to
the second summand.
We just notice that taking the residue for $t=0$
(that is, the coefficient of $t^{-1}$)
is the analog of the codimension-one condition $\d-\d' + |\b'-\b| = d-1$.

\

\ni
{\bf 3.3.\  Generalizations.\  }
The recursive procedure that we just described can easily be applied to
Hirzebruch surfaces, by replacing
the line $L$ with the exceptional curve $E$. Then one
obtains
a completely similar formula for curves of
any genus (see [V]). A more subtle problem is generalizing the method to higher
dimensional projective spaces.
More precisely, one could ask for the number of (smooth) curves of given
degree and genus in $\P ^n$
that pass
through the appropriate number of points, or
 more generally   which satisfy a given set of
intersection conditions with respect to linear subspaces  of varying dimension.
For example, it is not hard to show that,
in $\P ^3$,  there are finetely many rational curves of degree $d$ satisfying
$4d$  linear conditions (ie, passing through $4d$ points, or meeting $4d$
lines, and so on).

The natural way of generalizing our technique is to place the ``linear
conditions" on a fixed
hyperplane one at the time so as to get a recursion.
This has been shown to work  by Vakil, for curves of genus $0$ in $\P ^n$.
The higher genus case is
 still under investigation.

\

\

\cl{4. EXAMPLES}

\

\

Consider the following well known

\proclaim Proposition.
There are 12 rational plane cubics passing through 8 general points.

\ni
In symbols: $N_0(3) = N^{3,1}(0,3) = 12$.
We give here a few different ways to prove such a result, to illustrate how
the  techniques we talked about
work. We first prove it by classical arguments (two of them) in 4.1, then
we use the cross ratio method and
the rational fibration method in 4.2 (these are ad-hoc techniques for
rational curves, ispired by the First
Reconstruction Theorem of [KM]). Finally, we use the recursive procedure of
[CH3] which led to
the Theorem stated in the previous chapter.

We let $V=V_0(3)\subset \P ^9$ and $N=N_0(3)$; we have $\dim V = 8$.

\ni
Remark. The following fact about the geometry of $V$ will be needed (cf.
[CH3] and loc. cit.). $V$ is
irreducible and smooth at its general point (corresponding to an
irreducible  cubic with one node). $V$
contains a codimension 1 irreducible subvariety $W$ whose general points
parametrize reducible cubics given by
the union of a conic and a line. This is the unique degenerate locus of $V$
having codimension 1 and
parametrizing reducible curves.
$V$ is singular along $W$, looking like two smooth sheets crossing
transversally. The two sheets correspond to
deformations of the reducible curves that are locally trivial at either one
of the nodes.

\

\ni
{\bf 4.1.\  Classical proofs.\  }
These proofs are  well known.
We  actually prove the more general fact that the
degree of the Severi variety of curves of degree $d$ with one node is
$3(d-1)^2$.
So, let $V=V^{d,1}[0,d]$  and $N$ be its degree. $V$ has codimension $1$ in
the space $\Pr$
of all curves of
degree
$d$, so that $N= \ell \cap V$ where $\ell$ is a general line in $\Pr$.

We can identify $\ell$ with a general pencil of curves of degree $d$; that
is, a family given by a polynomial
equation $F(t_0,t_1;x_0,x_1,x_2)=0$ where $F$ is homogeneous of degree $1$
in $t_0,t_1$ and of degree $d$ in
$x_0,x_1,x_2$. We obtain this way a family
$$
\Y \la \P ^1 \cong \ell
$$
where $\Y \subset \P ^1 _{t_0,t_1}\times \Pl _{x_0,x_1,x_2}$.
Then $N$ corresponds exactly to the total number of nodes of the fibers of
such a family.

There are now two ways of computing such a number, an algebraic way and a
topological way.
Algebraically, the nodes of the fibers are given by the (non-zero)
solutions of the system
$$
F_{x_0}=F_{x_1}=F_{x_2}=0,
$$
that is, by the number of points in which the three surfaces given by
$F_{x_i}=0$ in
$\P ^1 _{t_0,t_1}\times \Pl _{x_0,x_1,x_2}$ intersect. Now, these surfaces
are of class $h_1+(d-1)h_2$, where
$h_i$ is the pull-back of the generator of $\Pic (\P ^i)$, $i=1,2$. Hence
the  basic relations
$h_1^3-h_1^2h_2=h_2^3=0$ and $h_1h_2^2=1$ imply that
$(h_1+(d-1)h_2)^3=3(d-1)^2$, which is what we wanted to prove. We obtain
the proposition as a special case.

The same result can be obtained by computing the Euler characteristic of
$\Y$ in two different ways.
First, $\Y$ is the blow up of $\Pl$ at the $d^2$ base points of the pencil,
hence
$$
\chi (\Y ) = 3+ d^2.
$$
On the other hand, the family $\Y \la \P ^1$ has general fiber $F$ a
Riemann surface of
topological
characteristic
$\chi (F) = 2-2g=2-2{d-1 \choose 2} $. If $N$ is the total number of nodes,
we have (cf. [GH] Chapter 4)
$$
\chi (\Y )=\chi (F)\chi (\P ^1)+N
$$
which implies $N=3(d-1)^2$.

\

\ni {\bf 4.2.\ The cross ratio and the rational fibration methods.\  }
Here we prove the proposition by using the fact that the curves in question
are rational
(see also [DI]).
Fix $7$ general points in the plane and let $\Gamma \subset V$ be the
irreducible curve parametrizing all
 nodal
cubics through such points. Let $\X \la \Gamma$ be the corresponding
family, so that $ \X \subset \Gamma
\times
\Pl$. Observe that the family has exactly ${7 \choose 2}$ reducible fibers,
corresponding to all
reducible cubics of type $C_1\cup C_2$, where $C_1$ is a line through two
of the base points and $C_2$ is the
conic through the remaining five points.
If $t\in \Gamma $ is a point such that the fiber $X_t$ is one of these
reducible curves, then $t$ will be a
node of $\Gamma$ by the introductory remark. We then let $B$ be the
normalization of $ \Gamma$ and we let
$\Y$ be the normalization of the base change family $\X \times _{\Gamma }B$.
It turns out that $\Y$ is  a smooth surface and that all of the fibers of
$\Y \la B$ are at most nodal.
We do not really need this; if it were not the case that all fibers were
nodal (there may be a priori some
cuspidal curves) we could make a base change and perform semistable
reduction; if $\Y$ were not smooth, we
could take its minimal desingularization. None of these two operations would affect the rest of the procedure.

Finally  $\Y \la B$ is a  family of generically smooth rational curves,
having $2{7 \choose 2}$ reducible
nodal fibers made of two (smooth, rational) components, and no other
singular fiber.
We let $\pi :\Y \la \Pl$ be the natural map.

This part of the set-up is common to both the cross ratio method and the
rational fibration method.
Now we concentrate on the first. From
$\Y \to B$ we want to obtain a family of rational curves having $4$
sections. We  do that as follows:
the first two sections will be two of the $7$ base points, call them
 $p_1$ and $p_2$; the other two sections
will be given by intersecting the curves of the family with $2$ fixed lines
$L_3$ and $L_4$ (which are
chosen to be general with respect to the base points of the family). This
will clearly be possible after a
base change of order $9$; in fact, $L_3$ intersects each curve of the
family in $3$ points, hence $\pi ^*L_3$
is a curve in
$\Y$ wich is a $3$ to $1$ cover of $B$. Therefore  if we perform the  base
change of degree $3$
given by
$\pi ^*L_3 \to B$ we can  define a (single-valued) section of $\Y \times _B
\pi ^*L_3 $ representing the
intersection of $L_3$ with the curves of the family. Then we repeat the
same process to obtain a
(single-valued) section out of the intersection with $L_4$. In total, we
made a degree $9$ base change
$A\to B$ and we get a new family $\W \to A$ where $\W := \Y \times _B A$.
We call again $\pi $ the natural
map from $\W$ to $\Pl$.
By construction, this is a family of generically smooth rational curves,
all of whose
fibers are at most nodal, and such that there are
$4$ sections $p_i:A \to \W$, for $i = 1,2,3,4$.  These sections are such
that the first two correspond to
$p_1$ and $p_2$ respectively, while $\pi (a) \in L_i$ for $i=3,4$ and for
all $a\in A$.
The reader familiar with the theory of moduli of curves
will immediately see that this family has a canonical
morphism
$\phi :A\la \overline{M_{0,4}}$ where $\overline{M_{0,4}}$  denotes the
Deligne-Mumford moduli space
of stable curves of genus $0$ with $4$ marked points.
There is a completely equivalent way of describing this map. Define
$\varphi :A\la \P ^1$ via the
cross ratio:
$$
\varphi (a) = {\bigl( p_1(a)-p_3(a)\bigr) \bigl( p_2(a)-p_4(a)\bigr) \over
\bigl( p_1(a)-p_2(a)\bigr)
\bigl( p_3(a)-p_4(a)\bigr) }.
$$
Notice that $\overline{M_{0,4}}\cong \P ^1$ so that if we want $\phi$ to
coincide with $\varphi$ we just have
to identify the $3$ boundary points of $\overline{M_{0,4}}$ with the
degenerate values $0$, $1$, and $\infty$
of the cross ratio as follows. The value $0$  is identified with the point
corresponding to the
isomorphism class of stable, nodal, rational  curves
 having $p_1$ and $p_3$ on one component and $p_2$, $p_4$ on the other
component (we call these
curves of type
$(13,24)$).  The value $\infty$   with the stable nodal rational  curve
 having $p_1$ and $p_2$ on one component and $p_3$, $p_4$ on the other
component.

Now we  prove the Proposition using the basic
fact that $\deg \varphi ^*(0)=\deg \varphi ^*(\infty )$.

To compute the degree of the zero divisor of $\varphi$, we know that the
value $0$ is achieved on curves
of type $(13,24)$
(notice  that, by the genericity assumption, the sections $p_1 $ and $p_3$
are disjoint, and so are the
sections
$p_2 $ and $p_4$). To
count them, we observe that there are
$10$ reducible cubics through the seven base points, such that $p_1$ and
$p_2$ lie on different components: $5
$ have $p_1$ on the line  (because we need one more point, out of the
remaining $5$ base points, to pin down
the line) and
$5$ have
$p_1$ on the conic (as before the line is determined by $p_2$ and any one
of the remaining $5$ base points).
The corresponding
$10$ points of
$\Gamma$ will be nodes, by the remark at the beginning of the chapter,
correspondingly we get $20$ points on
$B$. Then we have to take into account the base change of order $9$, and,
more important, the fact that $L_i$
meets the conic in two points.
Finally we get
$$
\deg \varphi ^*(0) = 9\cdot 2\cdot 20.
$$

The poles of $\varphi$ are assumed  on points of $A$ corresponding to
reducible fibers of type
$(12,34)$ and
also on points of $A$ where the section $p_3$ crosses the section $p_4$
(the sections $p_1$and $p_2$ do not intersect). Notice that $p_3 $ and
$p_4$ intersect in all points $a\in
A$ such that the fiber $W_a$ is mapped by $\pi$ to a nodal cubics passing
through the $7$ base points and
through the point of intersection of $L_3$ and $L_4$. There are $N$ such
plane curves, and hence $9N$
corresponding points in $A$, where $\varphi $ has a pole.
To count the fibers of type $(12,34)$ we proceed as before: there is a unique
 plane reducible cubic of our family having $p_1$ and $p_2$  on a line,
then we have to account for the
node of $\Gamma$, for the base change of order $9$, and for the fact that
$L_3$ and $L_4$ meet the conic in
two points. We get a contribution of $2\cdot 2\cdot 2 \cdot 9$ to $\deg
\varphi ^*(\infty )$.
Finally, we find ${5 \choose 3}$ curves having $p_1$ and $p_2$ on a conic,
and a total contribution of
${5 \choose 3}\cdot 2 \cdot 9$  to $\deg  \varphi ^*(\infty )$.
We get
$$
\deg  \varphi ^*(\infty ) = 9N+2\cdot 2\cdot 2 \cdot 9+{5 \choose 3}\cdot 2
\cdot 9=
\deg \varphi ^*(0) =9\cdot 2\cdot 20
$$
which yields $N=12$.

The rational fibration method is somewhat simpler, because it does not
involve any further construction, once
we arrive at the above family $\Y  \la B$.
We let $Y$ be the class of the fiber, so that $Y^2=0$. Then we let $A$ be
the class of a section
corresponding to one of the seven base points, call it $q$. Then $(A\cdot Y)=1$.
Let $B'\subset B$ be the set of points such that the corresponding fiber is
reducible
(so that $B'$ contains exactly $2\cdot {7 \choose 2}$ points), for $b\in B'$ let
$Z_b$ and $W_b$ be the two components. Choose the names $Z$ and $W$ so that
$(A\cdot Z_b)=1$ and
$(A\cdot W_b)=0$ (in other words, $A$ always intersects the $Z_b$
component, for every $b\in B'$).
Then the classes
$Y$, $ A$, $\{ W_b \}_{b\in B'}$
generate the N\'eron-Severi group of $Y$, which is to say that we can
compute intersection numbers by
expressing every other class as a combination of these generators.
Notice that we have $W_b^2=-1$ and  $(A\cdot W_b) = (Y\cdot W_b)=0$.

Let $L$ be the hyperplane class in $\Pl$.
Then we  prove the proposition using
$$
(\pi ^*L\cdot \pi ^*L)=(\deg \pi )(L\cdot L)=N.
$$
Let us write
$$
\pi ^*L=c_YY+c_AA+\sum _{b\in B'}c_bW_b
$$
and now notice that

$$
\eqalign{&(i)   \quad  (\pi ^*L\cdot Y) = 3  \quad \Longrightarrow \quad
c_A=3 \cr
&   (ii)  \quad  (\pi ^*L\cdot A) = 0  \quad  \Longrightarrow \quad
c_Y=-3A^2 \cr
& (iii)   \quad (\pi ^*L\cdot W_b) = (L\cdot \pi _*W_b) = \deg \pi _*W_b
\quad \Longrightarrow
\quad c_b=-\deg \pi _*W_b.\cr}
$$
Using all of these relations we obtain
$$
(\pi ^*L\cdot \pi ^*L)=-9A^2-\sum _{b\in B'}(\deg \pi _*W_b)^2.
$$
To compute $A^2$, pick another of the base points, call it $q'$, and let
$A'$ be the corresponding section of
the family.
Since $A^2 = (A')^2$ and $(A\cdot A')=0$ we get $2A^2= (A-A')^2$. Now
$A-A'$ is supported  exactly on those
reducible fibers such that $q$ and $q'$ lie on different components (in
other words, $q'\in Z_b$). The number
of them is $2\cdot 10$, in fact there will be $5$ curves of the family such
that $q$ lies on a line and $q'$
on a conic, and $5$ such that $q$ lies on a conic and $q'$ on a line; the
factor of $2$ comes from the fact
that $\Gamma$ has a node at such curves.
So we conclude that $A^2=10$.

The last thing to compute is $\sum (\deg \pi _*W_b)^2$. This amounts to
count how many reducible fibers have
$q$ on the line, and how many have $q$ on the conic. Clearly $q$ is on a
line for $2\cdot 6$ fibers and on a
conic for $2\cdot {6 \choose 2}$ fibers.
Hence $\sum (\deg \pi _*W_b)^2=2\cdot 6\cdot 4 + 2\cdot {6 \choose 2}=78$
and we conclude
$N=9\cdot 10 - 78 =12$.

\

\ni {\bf 4.3.\  The proof using the Theorem of [CH3].\  }
As we said,  for the inductive technique we fix a line $L$ in $\Pl$ and we
pick points $p_1, p_2, p_3,
p_4,....$ on
$L$; then we study successive scheme-thoretic intersections
$$V_t=V\cap (\cap _{i=1}^tH_{p_i})$$
for which, of course, $\deg V= \deg V_t$.

By dimension count we have $V\cap H_{p_1} = V^{3,1}[(1),(2)]$ and
$V_2=V\cap H_{p_1}\cap H_{p_2} = V^{3,1}[(2),(1)]$ and these intersections
are transverse, hence there are no
extra factors in the degree computation.
The next step gives two components of dimension $5$, namely
$$
V_3 = V^{3,1}[(3),(0)]\cup W
$$
where the general point of $W$ is a reducible cubic $X=L\cup C$ where $C$
is a conic.
Hence $\deg W =1$, but there will be a coefficient of $2$ in the degree
computation.
This is because $V_2$ is singular along  $W$, by the remark at the
beginning of this chapter. This coefficient
$2$ corresponds to ${\b ' \choose \b}$ in the formula.

Now we have to compute the degree of
$V' := V^{3,1}[(3),(0)]$. By intersecting
with a fourth hyperplane we now get five different irreducible components of
dimension $4$, three of which are of the same type. Let us write
$$
V'\cap H_{p_4} = U_1\cup U_2 \cup U_3 \cup W' \cup W''.
$$
The general point of $U_i$ is a curve $X=L\cup C$ where $C$ is a conic through
$p_i$. To understand this degeneration in more detail, consider a generic
one-parameter family of curves of
$V'$ degenerating to $X$. This is a family whose general fiber has one node
while the special fiber
has two, one of which is $p_i$.

We claim that the limit on $X$ of the node of the general fiber must be $p_i$.
To see this, we normalize the total space of the family so as to obtain a
family of generically smooth
curves. Let $Y$ be the curve lying over $X$, so that $Y$ is the partial
normalization at the limit node.
Now the total space of this normalized family has an obvious map to $\Pl$,
call it $\pi$.
The preimage $\pi ^{-1}L$ of the fixed line $L$ cannot have isolated
points, because $L$ is irreducible of
codimension $1$. Hence  $\pi ^{-1}L=L'\cup S_1\cup S_2\cup S_3$  where $L'$
is a subcurve of $Y$ and
$S_i$ is a curve such that $\pi S_i=p_i$.
Let $q\in Y$ be the point lying over the node of $X$ such that $q\not\in
L'$. Clearly $q\in \pi ^{-1}L$,
therefore there must be a curve of $\pi ^{-1}L$ passing through $q$. Such a
curve can only be one of the
$S_i$, and this implies our claim.

We conclude that, in the degree computation, we shall find no further
coefficient, as the intersection is
transverse and $V'$ is smooth along $U_i$.
We have of course
$\deg U_i = N^{2,0}[(1),(1)]=1$.

The general curve  parametrized by $W'$ is the union of $L$ and a conic $C$
tangent to $L$ at some unassigned point. The degree of the variety of
conics tangent to $L$ is
$2$, and $H_{p_4}$ is tangent to $V'$ along $W'$, hence we get a factor of
$2$ which
correspond to
$\N ^{\b ' - \b}$ in the Theorem.

Finally $W''$ parametrizes curves of type $X=L\cup L_1 \cup L_2$ where
$L_i$ is a line; the same argument
used before shows that the limit node of such a curve can only be $L_1\cap
L_2$, so that $V'$ is smooth
along $W''$ and the intersection is transverse.
 The degree of the variety of nodal conics is $3$.

We can summarize the above description of the (scheme intersection)
$V^{3,1}[(3)(0)] \cap H_{p_4}$ by looking at the following three
possibilities for the limit of the node of
the general fiber. If  such a limit is an assigned point, say $p_i$, then
the limit curve lies in $U_i$; if
the limit node is an unassigned point of $L$, then the limit curve
is in $W'$; finally, if the limit node is
not on $L$, we get a curve in $W''$.

We conclude that $N^{3,1}[(3)(0)]=3+2\cdot 2 +3 = 10$ and hence $N_0(3)=12$.
The reader might find the following explicit formula  useful
$$
N_0(3)=2\cdot N^{2,0}[(0),(2)]  + 3\cdot N^{2,0}[(1),(1)] + 2\cdot
N^{2,0}[(0),(0,1)]
+N^{2,1}[(0),(0)]
$$
where the coefficient $3$ in front of the second summand corresponds to
${\a \choose
\a '}$.

\

\

\cl{REFERENCES}

\

\

\ni [CH1] L. Caporaso and  J. Harris, {\it Parameter spaces for curves on
surfaces and
enumeration of rational curves.} \ (1995) Preprint,  alg-geom 9608023.

\ni [CH2] L. Caporaso and  J. Harris, {\it Enumeration of rational curves: the
rational fibration method.} \ (1995) Preprint,  alg-geom 9608024.

\ni [CH3] L. Caporaso and  J. Harris, {\it Counting plane curves of any genus}
(1996) Preprint,  alg-geom 9608025.

\ni [Ch] Y. Choi, {\it Severi degrees in cogenus 4} (1996) Preprint,
alg-geom 9601013.

\ni [CM] B. Crauder  and  R. Miranda, {\it Quantum cohomology of rational
surfaces.}
The Moduli Space of Curves, Progress in mathematics
129,   Birkhauser (1995)  pp. 33-80.

\ni [DI]  P. Di Francesco and C. Itzykson,
{\it  Quantum intersection rings.} The Moduli Space of Curves, Progress in
mathematics
129,   Birkhauser (1995)  pp. 81-148.

\ni
[FP] W. Fulton and R. Pandharipande, {\it Notes on stable maps and quantum
cohomology.} (1995) Preprint.

\ni [G] E. Getzler, {\it Intersection theory on $\overline{M}_{1,4}$ and
Gromov-Witten invariants in genus 1.}
MPI Preprint (1996)

\ni [H] \  J. Harris, {\it On the Severi problem}.  \ Invent. Math. 84
(1986), pp. 445-461.

\ni [KP] S. Kleiman and R. Piene,
Private communication.

\ni [Ko] \ J. Koll\'ar, {\it Rational Curves on Algebraic Varieties},
Springer 1996.

\ni [K] M. Kontsevich, {\it Enumeration of rational curves via torus action}.
The Moduli Space of Curves, Progress in mathematics 129,   Birkhauser
(1995)  pp.
335-368.

\ni [KM] M. Kontsevich and Y. Manin,  \ {\it Gromov-Witten classes, quantum
cohomology
and enumerative geometry.} Commun.Math.Phys. 164 (1994), pp. 525-562 and
hep-th/9402147.

\ni [P] R.  Pandharipande,
{\it Intersections of $\Q$-divisors on Kontsevich's moduli space

\ni
$\overline{M}_{0,n}(\P^r, d)$ and enumerative geometry}. \ 1995,
Preprint.

\ni [R1] Z. Ran, {\it Enumerative geometry of singular plane curves.}
  Invent. Math. 97 (1989), pp. 447-465.

\ni [R2] Z. Ran, {\it On the  quantum cohomology of the plane old and new.}
(1995) Preprint,  alg-geom 9508011.

\ni [RT] Y. Ruan and G. Tian, {\it A mathematical theory of Quantum Cohomology}
  J. Diff. Geom. 42 No. 2 (1995) pp. 259-367.

\ni [S] F. Severi, {\it Vorlesungen uber algebraische Geometrie} Anhang F.
Leipzig: Teubner 1921.

\ni [V] R. Vakil, {\it Curves on rational ruled surfaces} (1996) Preprint.

\ni [Va] I. Vainsencher,  {\it Counting divisors with prescribed
singularities.} Trans AMS 267 (1981), 399-422.

\ni [W] E. Witten, {\it Two dimensional gravity and intersection theory on
moduli space.} Surveys in Diff.
Geom. 1 (1991) pp. 243-310.

\

\noindent Mathematics department, Harvard University, 1 Oxford st.,
Cambridge MA 02138, USA

\noindent caporaso@abel.harvard.edu

\end